\documentclass[leqno]{article}
\usepackage{amsmath}
\usepackage{amssymb}
\def\C{{\mathbb{C}}}

\def\R{{\mathbb{R}}}
\def\N{{\mathbb{N}}}
\def\Z{{\mathbb{Z}}}
\def\P{{\mathbb{P}}}

\def\RE{\operatorname{Re}}
\def\ds{\displaystyle}

\newtheorem{theorem}{Theorem}[section]
\newtheorem{corollary}{Corollary}[section]
\newtheorem{example}{Example}[section]

\newtheorem{remark}{Remark}[section]

\begin{document}

\title{\bf On the binomial convolution of arithmetical
functions \footnote{Journal of Combinatorics and Number Theory,
Volume 1, Issue 1, 2009, pp. 31--48}}
\author{{\sc L\'aszl\'o T\'oth} \\ University of P\'ecs, Institute of Mathematics and Informatics \\
Ifj\'us\'ag u. 6, 7624 P\'ecs, Hungary \\
{E-mail address: ltoth@ttk.pte.hu} \\ \\
{\sc Pentti Haukkanen} \\ Department of Mathematics and Statistics,\\
FI-33014 University of Tampere, Finland \\ E-mail address:
pentti.haukkanen@uta.fi}
\date{ }
\maketitle

{\bf Abstract:} Let $n=\prod_p p^{\nu_p(n)}$ denote the canonical
factorization of $n\in \N$. The binomial convolution of arithmetical
functions $f$ and $g$ is defined as $(f\circ g)(n)=\sum_{d\mid n}
\left( \prod_p \binom{\nu_p(n)}{\nu_p(d)} \right) f(d)g(n/d),$ where
$\binom{a}{b}$ is the binomial coefficient. We provide properties of
the binomial convolution. We study the $\C$-algebra $({\cal
A},+,\circ,\C)$, characterizations of completely multiplicative
functions, Selberg multiplicative functions, exponential Dirichlet
series, exponential generating functions and a generalized binomial
convolution leading to various M\"obius-type inversion formulas.
Throughout the paper we compare our results with those of the
Dirichlet convolution $*$. Our main result is that $({\cal
A},+,\circ,\C)$ is isomorphic to $({\cal A},+,*,\C)$. We also obtain
a ``multiplicative'' version of the multinomial theorem.



\section{Introduction}

Let ${\cal A}$ denote the set of arithmetical functions $f:\N \to
\C$. It is well known that ${\cal A}$ is a $\C$-algebra under the
linear operations and the Dirichlet convolution defined by
\begin{equation}
(f*g)(n)=\sum_{d\mid n} f(d)g(n/d).   \label{Dirichlet convo}
\end{equation}
This $\C$-algebra is denoted as $({\cal A},+,*,\C)$. It is
isomorphic to the $\C$-algebra of formal Dirichlet series
$D(f,s)=\sum_{n=1}^{\infty} f(n)n^{-s}$, denoted as $({\cal
D},+,\cdot,\C)$. Then $D(f*g,s)=D(f,s)D(g,s)$ and the mapping
$f\mapsto D(f,s)$ serves as an isomorphism. Furthermore, $({\cal
A},+,*)$ is an integral domain; it is even a unique factorization
domain, cf. \cite{McC1986,Ten1995}.

Let $n=\prod_p p^{\nu_p(n)}$ denote the canonical factorization of
$n\in \N$. The binomial convolution of arithmetical functions $f$
and $g$ is defined as
\begin{equation}\label{binomial convo}
(f\circ g)(n)=\sum_{d\mid n} \left( \prod_p
\binom{\nu_p(n)}{\nu_p(d)} \right) f(d)g(n/d),
\end{equation}
where $\binom{a}{b}$ is the binomial coefficient. This convolution
appears in the book by {P. J. McCarthy} \cite[p. 168]{McC1986}, and
its basic properties were investigated by {P. Haukkanen}
\cite{Hau1996}, see also \cite[p. 116]{SC2004}. It was pointed out
in \cite{Hau1996} that the binomial convolution \eqref{binomial
convo} possesses properties analogous to those of the Dirichlet
convolution \eqref{Dirichlet convo}. For example, $(\{f\in {\cal A}:
f(1)\ne 0\},\circ)$ is a commutative group and the set of all
multiplicative arithmetical functions $f$ forms a subgroup of this
group. It is remarkable that the binomial convolution preserves
complete multiplicativity of arithmetical functions. The set of
completely multiplicative functions forms a subgroup of the group of
multiplicative functions under the binomial convolution. This is not
the case for the Dirichlet convolution. Note also that the inverse
of the function $I(n)=1$ ($n\in \N$) under the binomial convolution
is the Liouville function $\lambda(n)=(-1)^{\Omega(n)}$, where
$\Omega(n)=\sum_p \nu_p(n)$, while the inverse of $I$ under the
Dirichlet convolution is the M\"obius function $\mu$. The
arithmetical function $\delta$ defined as $\delta(1)=1$,
$\delta(n)=0$ for $n>1$ serves as the identity under both the
binomial and Dirichlet convolution.

In this paper we provide further properties of the binomial
convolution. We study the $\C$-algebra $({\cal A},+,\circ,\C)$,
characterizations of completely multiplicative functions, Selberg
multiplicative functions, exponential Dirichlet series, exponential
generating functions and a generalized binomial convolution leading
to various M\"obius-type inversion formulas. Throughout the paper we
compare our results with those of the Dirichlet convolution. Our
main result is that $({\cal A},+,\circ,\C)$ is isomorphic to $({\cal
A},+,*,\C)$. We also obtain a ``multiplicative'' version of the
multinomial theorem.


\section{The algebra $({\cal A},+,\circ,\C)$}

It is easy to see that $({\cal A},+,\circ,\C)$ is a $\C$-algebra. In
this section we show that $({\cal A},+,\circ,\C)$ is isomorphic to
$({\cal A},+,*,\C)$ and compare the expressions and inverses of the
convolutions in these algebras.

The function $\xi$ defined by $\xi(n)=\prod_p \nu_p(n)!$ plays a
crucial role in connections between the Dirichlet and binomial
convolution. We recall that an arithmetical function $f$ is said to
be multiplicative if $f(1)=1$ and $f(mn)=f(m)f(n)$, whenever $(m,
n)=1$, and completely multiplicative if $f(1)=1$ and
$f(mn)=f(m)f(n)$ for all $m$ and $n$. The function $\xi$ is
multiplicative and prime independent. Also, $\xi(n)=1$ if and only
if $n$ is squarefree. The function $\xi$ is not completely
multiplicative. However, $\xi(m)\xi(n)\mid \xi(mn)$ for all $m,n\ge
1$. In fact, one has $a!\, b! \mid (a+b)!$ for all $a,b\ge 1$ and
hence $\xi(p^a)\xi(p^b)=a!\, b! \mid (a+b)! = \xi(p^{a+b})$ for all
prime powers $p^a$, $p^b$. Therefore $\xi(m)\xi(n)\mid \xi(mn)$ and
in particular $\xi(m)\xi(n)\le \xi(mn)$.

\begin{theorem}
The algebras $({\cal A},+,\circ,\C)$ and $({\cal A},+,*,\C)$ are
isomorphic under the mapping $\ds f\mapsto\frac{f}{\xi}$.
\end{theorem}

{\bf Proof.} It is easy to see that the mapping $\ds
f\mapsto\frac{f}{\xi}$ is a bijection on ${\cal A}$. Furthermore,
$\ds f+g\mapsto \frac{f+g}{\xi}= \frac{f}{\xi}+\frac{g}{\xi}$ and
$\ds sf\mapsto \frac{sf}{\xi}=s\frac{f}{\xi}$. Moreover,
\begin{equation*}
(f\circ g)(n)=\sum_{de= n} \left(\prod_p \frac{\nu_p(n)!}{\nu_p(d)!
\, \nu_p(e)!} \right) f(d)g(e) =\xi(n) \sum_{de= n}
\frac{f(d)}{\xi(d)}\cdot \frac{g(e)}{\xi(e)} =\xi(n)
\left(\frac{f}{\xi}*\frac{g}{\xi}\right)(n),
\end{equation*}
showing that
\begin{equation}
f\circ g =\xi \left(\frac{f}{\xi}*\frac{g}{\xi}\right). \label{isom}
\end{equation}

Thus, $\ds f\circ g \mapsto \frac{f\circ g}{\xi}=
\frac{f}{\xi}*\frac{g}{\xi}$. This shows that $\ds
f\mapsto\frac{f}{\xi}$ is an algebra isomorphism.

\begin{corollary}
$({\cal A},+,\circ)$ is an integral domain. It is even a unique
factorization domain.
\end{corollary}

Equation \eqref{isom} expresses the binomial convolution in terms of
the Dirichlet convolution. On the other hand, we have
\begin{equation}\label{disom}
f*g =\frac{f\xi\circ g\xi}{\xi}
\end{equation}
or $ (f*g)\xi =f\xi \circ g\xi $ for all $f,g\in {\cal A}$.

We next write the binomial inverse in terms of the Dirichlet inverse
and vice versa. Let $f^{-1\circ}$ and $f^{-1*}$ denote the inverses
of $f$ under the binomial convolution $\circ$ and the Dirichlet
convolution $*$, respectively. They exist if and only if $f(1)\ne
0$.

\begin{theorem} \label{th:inv}
For any $f\in {\cal A}$ with $f(1)\ne 0$,
\begin{equation} \label{BinInv}
\ds f^{-1\circ}=\xi\left(\frac{f}{\xi}\right)^{-1*}
\end{equation}
and
\begin{equation} \label{DirInv}
\ds f^{-1*}=\frac{(\xi f)^{-1\circ}}{\xi}.
\end{equation}
\end{theorem}

{\bf Proof.} We have $\ds f\circ f^{-1\circ}=\delta$ and from
\eqref{isom}  we obtain $\ds
\xi\left(\frac{f}{\xi}*\frac{f^{-1\circ}}{\xi}\right)=\delta$ or
$\ds \frac{f}{\xi}*\frac{f^{-1\circ}}{\xi}=\delta$. This proves
Theorem \ref{th:inv}

\begin{example}\upshape
For $f=\xi$ we have $\xi^{-1\circ}=\xi I^{-1*}=\xi \mu =\mu$. Thus,
the inverse of $\xi$ under the binomial convolution is the M\"obius
function $\mu$, see \cite[p. 213]{Hau1996}.
\end{example}

A further result involving the binomial and Dirichlet inverses is
presented below. This result involves multiplicative functions.

\begin{theorem}\label{th:mult}
If $f$ is multiplicative and $f(p^a)=0$ for all prime powers $p^a$
with $a\ge 2$, then for every $n\ge 1$,
\begin{equation} \label{primepowerB}
f^{-1\circ}(n)=(-1)^{\Omega(n)}\xi(n)\prod_p f(p)^{\nu_p(n)}
=\lambda(n)\xi(n)\prod_p f(p)^{\nu_p(n)}
\end{equation}
and
\begin{equation} \label{primepowerD}
f^{-1*}(n)=(-1)^{\Omega(n)}\prod_p f(p)^{\nu_p(n)}
=\lambda(n)\prod_p f(p)^{\nu_p(n)}.
\end{equation}
\end{theorem}

{\bf Proof.} Let $p^a$ be a prime power with $a\ge 1$. Then
$0=(f\ast f^{-1\ast})(p^a)=f^{-1\ast}(p^a) +
f^{-1\ast}(p^{a-1})f(p)$, and thus
$f^{-1\ast}(p^a)=-f^{-1\ast}(p^{a-1})f(p)$. This recursion gives
$f^{-1\ast}(p^a)=(-1)^a\, f(p)^{a}$. Thus \eqref{primepowerD} holds
for all prime powers and therefore by multiplicativity it holds for
all positive integers. Equation \eqref{primepowerB} follows from
\eqref{BinInv} and \eqref{primepowerD}.

\begin{example}\upshape
For $f=\mu^2$ we have $f(p)=1$ for all primes $p$ and
$f^{-1\circ}(n)=(-1)^{\Omega(n)}\xi(n)=\lambda(n)\xi(n)$. For
$f=\mu$ we have $\mu^{-1\circ}(n)=(-1)^{\Omega(n)}\xi(n)\prod_p
(-1)^{\nu_p(n)}= (-1)^{\Omega(n)}\xi(n)(-1)^{\Omega(n)}=\xi(n)$.
This follows also from the result $\xi^{-1\circ}=\mu$. If $f(p)=r$
for all primes $p$, then
$f^{-1\circ}(n)=(-1)^{\Omega(n)}\xi(n)\prod_p r^{\nu_p(n)}=
(-r)^{\Omega(n)}\xi(n)$.
\end{example}

\begin{remark}\upshape
The function $f^{-1*}$ in \eqref{primepowerD} is completely
multiplicative for all $f$ satisfying the conditions in Theorem
\ref{th:mult}
\end{remark}

Equation \eqref{isom} can be extended to several functions. In fact,
from \eqref{isom} we obtain that $f\circ g \circ h = \xi
\left(\frac{f}{\xi}*\frac{g}{\xi}\right) \circ h = \xi
\left(\frac{f}{\xi}*\frac{g}{\xi}
* \frac{h}{\xi}\right)$ and in general for all $f_1,\ldots,f_k\in {\cal
A}$,
\begin{equation}
f_1\circ \cdots \circ f_k = \xi \left(\frac{f_1}{\xi}* \cdots *
\frac{f_k}{\xi}\right). \label{kfunct}
\end{equation}

This means that for every $n\in \N$,
\begin{equation}
(f_1\circ \cdots \circ f_k)(n)= \sum_{d_1\cdots d_k=n} \left(
\prod_p \binom{\nu_p(n)}{\nu_p(d_1),\ldots,\nu_p(d_k)}\right)
f_1(d_1)\cdots f_k(d_k), \label{kfunctions}
\end{equation}
involving multinomial coefficients. If $f_1,\ldots,f_k\in {\cal A}$
are multiplicative functions, then $f_1\circ \cdots \circ f_k$ is
multiplicative and
\begin{equation}
(f_1\circ \cdots \circ f_k)(n)= \sum_{d_1\cdots d_k=n} \left(\prod_p
\binom{\nu_p(n)}{\nu_p(d_1),\ldots,\nu_p(d_k)} f_1(p^{\nu_p(d_1)})
\cdots f_k(p^{\nu_p(d_k)})\right). \label{kmultip}
\end{equation}

Equation \eqref{disom} can also be extended to several functions. We
do not need these details in this paper.


\section{Completely multiplicative functions}

In \cite{Hau1996} the second author provides properties of
completely multiplicative functions with respect to the binomial
convolution. In this section we provide further properties of this
kind. In fact, we derive two characterizations of completely
multiplicative functions in terms the binomial convolution (Theorems
3.2 and 3.3). A large number of
characterizations of completely multiplicative functions in terms
the Dirichlet convolution have been published in the literature, see
e.g. \cite{A1971,H2001,LPW2002}. In Section 5 we find
the exponential Dirichlet series of completely multiplicative
functions.

We begin this section by deriving ``multiplicative'' versions of the
multinomial and binomial theorems from \eqref{kmultip}.

Let $f_1,\ldots,f_k$ be completely multiplicative aritmetical
functions. Then from \eqref{kmultip} we obtain
\begin{equation}
(f_1\circ \cdots \circ f_k)(n)= \sum_{d_1\cdots d_k=n} \left(\prod_p
\binom{\nu_p(n)}{\nu_p(d_1),\ldots,\nu_p(d_k)} f_1(p)^{\nu_p(d_1)}
\cdots f_k(p)^{\nu_p(d_k)}\right). \label{kcomplmultip}
\end{equation}

On the other hand, $f_1\circ \cdots \circ f_k$ is also completely
multiplicative and thus
\begin{equation}
(f_1\circ \cdots \circ f_k)(n)= \prod_p (f_1\circ \cdots \circ
f_k)(p)^{\nu_p(n)} = \prod_p (f_1(p)+\cdots+f_k(p))^{\nu_p(n)}.
\label{k}
\end{equation}

Now, suppose that $f_1,\ldots,f_k$ are prime independent completely
multiplicative functions,  that is, $f_1(p)=x_1,\ldots, f_r(p)=x_r$
for any prime $p$, where $x_1,\ldots,x_r$ are given complex numbers.
Then by \eqref{kcomplmultip},
\begin{equation}
(f_1\circ \cdots \circ f_k)(n)= \sum_{d_1\cdots d_k=n} \left(\prod_p
\binom{\nu_p(n)}{\nu_p(d_1),\ldots,\nu_p(d_k)} \right)
x_1^{\Omega(d_1)} \cdots x_k^{\Omega(d_k)}, \label{kpi-complmultip}
\end{equation}
and by \eqref{k}
\begin{equation}
(f_1\circ \cdots \circ f_k)(n)= (x_1+\cdots +x_k)^{\Omega(n)}.
\label{kpi-complmultip2}
\end{equation}

From \eqref{kpi-complmultip} and \eqref{kpi-complmultip2} we obtain
the following {\sl ``multiplicative'' version of the multinomial
theorem}. It reduces to the usual multinomial theorem if  $n$ is a
prime power.

\begin{theorem}\label{th:multmult}
For all complex numbers $x_1,\ldots,x_r$ and positive integers $n$,
\begin{equation}
\sum_{d_1\cdots d_k=n} \left( \prod_p
\binom{\nu_p(n)}{\nu_p(d_1),\ldots,\nu_p(d_k)} \right)
x_1^{\Omega(d_1)} \cdots x_k^{\Omega(d_k)} =(x_1+\cdots
+x_k)^{\Omega(n)}.  \label{multip-multinom}
\end{equation}
\end{theorem}

\vskip1mm For $k=2$ we obtain the following {\sl ``multiplicative''
version of the binomial theorem}.

\begin{corollary}\label{co:multbin}
For all complex numbers $x$ and $y$ and positive integers $n$,
\begin{equation}
\sum_{d\mid n} \left( \prod_p \binom{\nu_p(n)}{\nu_p(d)}\right)
x^{\Omega(d)} y^{\Omega(n/d)} =(x+y)^{\Omega(n)}.
\label{multip-binom}
\end{equation}
\end{corollary}

And, for $x_1=\cdots=x_k=1$, Theorem 3.1 becomes

\begin{corollary}
For all positive integers $k$ and $n$,
\begin{equation}
\sum_{d_1\cdots d_k=n} \prod_p
\binom{\nu_p(n)}{\nu_p(d_1),\ldots,\nu_p(d_k)} =k^{\Omega(n)}.
\label{multip-sum}
\end{equation}
\end{corollary}

We next write a characterization of completely multiplicative
functions in terms of binomial powers. A similar result in terms of
Dirichlet powers is presented in \cite{LPW2002}.

For $k\in \Z, k\ne 0$ let $f^{k\circ}$ denote the $k$-th power of
$f\in {\cal A}$ under the binomial convolution, i.e.
$f^{k\circ}=f\circ \cdots \circ f$ ($k$ times),
$f^{-k\circ}=f^{-1\circ}\circ \cdots \circ f^{-1\circ}$ ($k$ times)
for all $k>0$.

If $f$ is completely multiplicative, then by \eqref{k}
$f^{k\circ}(n)=k^{\Omega(n)}f(n)$ for all $k>0$ and $n\ge 1$. Also,
$f^{-1\circ}(n)= \lambda(n)f(n)$, therefore
$f^{-k\circ}(n)=k^{\Omega(n)} \lambda(n) f(n)=(-k)^{\Omega(n)}
f(n)$, that is $f^{k\circ}(n)=k^{\Omega(n)}f(n)$ for all $k\in \Z,
k\ne 0$ and $n\ge 1$. Defining $f^{(0)}=\delta$ and
$0^{\Omega(1)}=0^0=1$ this holds also for $k=0$.

The following is a sufficient condition for a multiplicative
function to be completely multiplicative.

\begin{theorem}\label{th:charS}
Let $f$ be a multiplicative function. If there is an integer $k\in
\Z$, $|k|\ge 2$ such that $f^{k\circ}(n)=k^{\Omega(n)} f(n)$ for all
$n\ge 1$, then $f$ is completely multiplicative.
\end{theorem}

{\bf Proof.} Suppose that $k\ge 2$. For any prime power $n=p^\nu$ we
have
\begin{equation}
f^{k\circ}(p^\nu)= k^\nu f(p^\nu). \label{proof1}
\end{equation}

According to \eqref{kmultip},
\begin{eqnarray*}
f^{k\circ}(p^\nu) &=& \sum_{\nu_1+\cdots+\nu_k=\nu}
\binom{\nu}{\nu_1,\ldots,\nu_k}
f(p^{\nu_1})\cdots f(p^{\nu_k})\\
&=& k f(p^\nu) + \sum_{\substack{\nu_1+\cdots+\nu_k=\nu\\
    \nu_1,\ldots,\nu_k<\nu}} \binom{\nu}{\nu_1,\ldots,\nu_k}
    f(p^{\nu_1})\cdots f(p^{\nu_k}).
\end{eqnarray*}

We show by induction on $m$ that $f(p^m)=f(p)^m$. This is true for
$m=1$. Assume that it holds for any $m<\nu$. Then
\begin{eqnarray*}
f^{k\circ}(p^\nu)
&=& k f(p^\nu) + \sum_{\substack{\nu_1+\cdots+\nu_k=\nu \\
\nu_1,\ldots,\nu_k<\nu}} \binom{\nu}{\nu_1,\ldots,\nu_k}
f(p)^{\nu_1}\cdots f(p)^{\nu_k}\\
&=& k f(p^\nu) + (\underbrace{f(p)+\cdots+f(p)}_k)^\nu-kf(p)^\nu = k
f(p^\nu) + k^\nu f(p)^\nu-kf(p)^\nu,
\end{eqnarray*}
by the multinomial formula. We obtain
\begin{equation}
f^{k\circ}(p^\nu)=k f(p^\nu) + (k^\nu-k)f(p)^\nu. \label{proof2}
\end{equation}

By \eqref{proof1} and \eqref{proof2}, $k^\nu f(p^\nu) = k f(p^\nu) +
(k^\nu-k)f(p)^\nu$ or $(k^\nu-k) f(p^\nu) = (k^\nu-k)f(p)^\nu$,
where $k^\nu- k\ne 0$ ($k\ge 2$). Therefore $f(p^\nu)=f(p)^\nu$.

Now, suppose that $k$ is negative. This case is reduced to what is
already proved. Let $k=-j$, $j\ge 2$. Then for every $n\ge 1$,
$f^{-j\circ}(n)=(-j)^{\Omega(n)}f(n)$, which can be written as
$(f^{-1\circ})^{j\circ}(n)=j^{\Omega(n)}\lambda(n)f(n)$ or
$(f^{-1\circ})^{j\circ}(n)=j^{\Omega(n)}f^{-1\circ}(n)$.
Since $f^{-1\circ}$ is multiplicative, it follows from the first
part of the proof that $f^{-1\circ}$ is completely multiplicative
and we conclude that $f$ is completely multiplicative.

\begin{remark}\upshape
If $f$ is completely multiplicative, then $f^{-1\circ}=\lambda f$.
The converse is not true: If $f$ is multiplicative and
$f^{-1\circ}=\lambda f$, then $f$ need not be completely
multiplicative. In fact, $f$ is given by
\begin{eqnarray*}
&&f(1)=1,\\
&&f(p^{2n-1})\ \ {\rm may\ be\ arbitrary,}\ \ n=1, 2,\ldots,\\
&&f(p^{2n})=-\sum_{k=1}^{n-1} {2n\choose k} (-1)^k f(p^k)
f(p^{2n-k})
            -{1\over 2}{2n\choose n} (-1)^n f(p^n)^2,
            \ \ n=1, 2,\ldots,
\end{eqnarray*}
see \cite[p. 215]{Hau1996}. Note that it is well known that if $f$
is multiplicative and $f^{-1\ast}=\mu f$, then $f$ is completely
multiplicative.
\end{remark}

\begin{remark}\upshape
The function $\mu^{k*}(n)$ plays the role of the function
$k^{\Omega(n)}$ in the Dirichlet powers. It is easy to see that
$k^{\Omega(n)}=\lambda^{k\circ}(n)$. Note that the functions
$\lambda^{k\circ}(n)$ form an infinite cyclic subgroup of the group
of completely multiplicative functions under the binomial
convolution, while the functions $\mu^{k*}(n)$ form an infinite
cyclic subgroup of the group of multiplicative functions under the
Dirichlet convolution, see \cite{BHWS2000}.
\end{remark}

It is well known that distributivity over the Dirichlet convolution
is a characterization of completely multiplicative functions, for
details, see e.g. \cite{H2001}. Similar results can also be derived
for the binomial convolution. As an example of such
characterizations we present the following basic result.

\begin{theorem}\label{th:chardistr}
Let $f$ be a multiplicative function. Then $f$ is completely
multiplicative if and only if $f(g\circ h)=fg\circ fh$ for all
$g,h\in {\cal A}$.
\end{theorem}

{\bf Proof.} The ``$\Rightarrow$'' direction is immediate.

We prove the ``$\Leftarrow$''  direction. Let $g(n)=\mu(n)$,
$h(n)=\xi(n)$. Then $\delta=f\delta = f(\mu\circ\xi)=f\mu \circ
f\xi$; hence $(f\mu)^{-1\circ}=f\xi$. On the other hand, $f\mu$ is
multiplicative and $(f\mu)(p^a)=0$ for all $a\ge 2$, and thus
according to Theorem 2.3,
$(f\mu)^{-1\circ}(n)=(-1)^{\Omega(n)}\xi(n)\prod_p
(-f(p))^{\nu_p(n)}= \xi(n) \prod_p f(p)^{\nu_p(n)}$.

We obtain $f(n)\xi(n)= \xi(n) \prod_p f(p)^{\nu_p(n)}$ or
$f(n)=\prod_p f(p)^{\nu_p(n)}$ for all $n\in\N$, showing that $f$ is
completely multiplicative.

\begin{remark}\upshape
For a construction which is similar to the binomial convolution of
completely multiplicative arithmetical functions see \cite[Section
4]{BNV2008}.
\end{remark}

%

\section{Selberg multiplicative functions}\label{Selberg}

An arithmetical function $F$ is said to be Selberg multiplicative if
for each prime $p$ there exists $f_p:\N_0\to\C$ with $f_p(0)=1$ for
all but finitely many $p$ such that
\begin{equation}\label{e:Sel-exp}
F(n)=\prod_{p} f_p(\nu_p(n))
\end{equation}
for all $n\in\N$. An arithmetical function $F$ is said to be
semimultiplicative if
$$
F(m)F(n)=F((m, n))F([m, n])
$$
for all $m, n\in\N$, where $(m, n)$ and $[m, n]$ stand for the gcd
and lcm of $m$ and $n$. It is known that  an arithmetical function
$F$ (not identically zero) is semimultiplicative if and only if
there exists a nonzero constant $c_F$, a positive integer $a_F$ and
a multiplicative function $F'$ such that
\begin{equation}\label{e:sem-rea}
F(n)=c_F F'(n/a_F)
\end{equation}
for all $n\in\N$. (We interpret that the arithmetical function $F'$
possesses the property that $F'(x)=0$ if $x$ is not a positive
integer.)  We will take $a_F$ as the smallest positive integer $k$
such that $F(k)\ne 0$. Note that $c_F=F(a_F)$. Furthermore, it is
known that an arithmetical function is Selberg multiplicative if and
only if it is semimultiplicative. Semimultiplicative functions $F$
with $F(1)\ne 0$ are known as quasimultiplicative functions.
Quasimultiplicative functions $F$ possess the property
$F(1)F(mn)=F(m)F(n)$ whenever $(m, n)=1$. Semimultiplicative
functions $F$ with $F(1)=1$ are the usual multiplicative functions.

A semimultiplicative function not identically zero possesses a
Selberg expansion \eqref{e:Sel-exp} as
$$
F(n)=F(a_F)\prod_{p}
\left(\frac{F(a_Fp^{\nu_p(n)-\nu_p(a_F)})}{F(a_F)}\right).
$$

A Selberg expansion \eqref{e:Sel-exp} of a multiplicative function
is
$$
F(n)=\prod_{p} F(p^{\nu_p(n)}).
$$

It is known that semimultiplicative functions form a commutative
semigroup with identity under the Dirichlet convolution and
\begin{equation}\label{eq:acFD}
a_{F\ast G}=a_Fa_G, \ c_{F\ast G}=c_Fc_G, \ (F\ast G)'=F'\ast G'.
\end{equation}

Quasimultiplicative functions form a commutative group under the
Dirichlet convolution. For material on Selberg multiplicative and
semimultiplicative functions we refer to
\cite{BHS,Hau,HHS,Rearick66a,Rearick66b,Selberg,S}.

We next prove that semimultiplicative functions form a commutative
semigroup with identity under the binomial convolution. In the proof
we use the following result. If $F$ and $G$ are semimultiplicative,
then $FG$ is also semimultiplicative. In particular, if $F$ is
semimultiplicative (not identically zero) and $f$ is multiplicative
with $f(a_F)\ne 0$, then
\begin{equation}\label{eq:acfF}
a_{fF}=a_F, \ c_{fF}=f(a_F)c_F, \ (fF)^\prime=\frac{f_{a_F}}{f(a_F)}
F^\prime,
\end{equation}
where $f_a(n)=f(an)$ for all $n\in \N$.

\begin{theorem}
Semimultiplicative functions form a commutative semigroup with
identity under the binomial convolution. Furthermore,
$$
a_{F\circ G}=a_Fa_G, \ c_{F\circ G}= c_Fc_G
\frac{\xi(a_Fa_G)}{\xi(a_F)\xi(a_G)}, \ (F\circ G)^\prime=
\frac{\xi_{a_Fa_G}}{\xi(a_Fa_G)\xi}
\left(\left(\frac{\xi(a_F)\xi}{\xi_{a_F}}F^\prime\right) \circ
\left( \frac{\xi(a_G)\xi}{\xi_{a_G}}G^\prime \right) \right),
$$
where $\xi_a(n)=\xi(an)$ for all $n\in \N$.
\end{theorem}

{\bf Proof.} Let $F$ and $G$ be semimultiplicative. We show that
$F\circ G$ is also semimultiplicative. We use the formula
\[
F\circ G =\xi\left(\frac{F}{\xi}*\frac{G}{\xi}\right).
\]

The function $\xi$ is a multiplicative function such that $\xi(n)\ne
0$ for all $n\in \N$. Therefore $F/\xi$ and $G/\xi$ are
semimultiplicative and thus $(F/\xi)* (G/\xi)$ and $\xi((F/\xi)*
(G/\xi))$ have the same property.

Now, since $\xi(n)\ne 0$ for all $n\in \N$, on the basis of
\eqref{eq:acFD} and \eqref{eq:acfF} we obtain
\[
a_{F\circ G}=a_{(F/\xi)*(G/\xi)}=a_{F/\xi}a_{G/\xi}=a_Fa_G.
\]

Furthermore,
\[
c_{F\circ G}=(F\circ G)(a_{F\circ G})= (F\circ G)(a_Fa_G)=
\xi(a_Fa_G)((F/\xi)*(G/\xi))(a_Fa_G),
\]
where, taking into account, that $a_F:=k$ and $a_G:=\ell$ are the
least numbers such that $F(k)\ne 0$ and $G(\ell)\ne 0$,
respectively, the second factor of the last expression is
\[
((F/\xi)*(G/\xi))(a_Fa_G)=\frac{F(a_F)}{\xi(a_F)}\cdot
\frac{G(a_G)}{\xi(a_G)}= \frac{c_Fc_G}{\xi(a_F)\xi(a_G)}.
\]

Finally, on the basis of \eqref{eq:acFD} and \eqref{eq:acfF} we have
\begin{eqnarray*}
(F\circ G)' &=&
\left(\xi\left(\frac{F}{\xi}*\frac{G}{\xi}\right)\right)'=
\frac{\xi_{a_Fa_G}}{\xi(a_Fa_G)}
\left(\frac{F}{\xi}*\frac{G}{\xi}\right)'\\
&=& \frac{\xi_{a_Fa_G}}{\xi(a_Fa_G)}
\left(\left(\frac{F}{\xi}\right)' *\left(\frac{G}{\xi}\right)'
\right)=
\frac{\xi_{a_Fa_G}}{\xi(a_Fa_G)}\left(\frac{\xi(a_F)}{\xi_{a_F}}F'
* \frac{\xi(a_G)}{\xi_{a_G}}G' \right)\\
&=& \frac{\xi_{a_Fa_G}}{\xi(a_Fa_G)\xi}
\left(\left(\frac{\xi(a_F)\xi}{\xi_{a_F}}F'\right) \circ \left(
\frac{\xi(a_G)\xi}{\xi_{a_G}}G' \right) \right).
\end{eqnarray*}

This completes the proof.

\medskip

Let $S$ denote the class of Selberg multiplicative functions $F$
such that there exists a universal $f:\N_0\to\C$ with $f(0)=1$ such
that
$$
F(n)=\prod_{p} f(\nu_p(n))
$$
for all $n\in\N$. The class $S$ is exactly the class of prime
independent multiplicative functions. It is known that $S$ forms a
subgroup of the commutative group of multiplicative functions under
the Dirichlet convolution and for $H=F\ast G$ we have
$h(r)=\sum_{i=0}^r f(i)g(r-i)$. In a similar way we can prove that
$S$ forms a subgroup of the commutative group of multiplicative
functions under the binomial convolution and for $H=F\circ G$ we
have $h(r)=\sum_{i=0}^r \binom{r}{i} f(i)g(r-i)$. Note that the
functions $\lambda^{k\circ}(n)$ form a subgroup of the group $(S,
\circ)$, while the functions $\mu^{k\ast}(n)$ form a subgroup of the
group $(S, \ast)$.

%

\section{Exponential Dirichlet series}\label{sec:EDS}

For an arithmetical function $f$ we define the (formal) exponential
Dirichlet series by
\[
\widetilde{D}(f,s)= D\left(\frac{f}{\xi},s\right)=
\sum_{n=1}^{\infty} \frac{f(n)}{\xi(n) n^s}.
\]

Then $\widetilde{D}(\xi,s)=\zeta(s)$ is the Riemann zeta function
(the  Dirichlet series of the constant function $1$) and we let
$\widetilde{D}(I,s)=\widetilde{\zeta}(s)$ denote the exponential
Dirichlet series of the constant function $1$. Exponential Dirichlet
series has not hitherto been investigated in the literature, while
the usual Dirichlet series is one of the most fundamental concepts
in analytic number theory, see e.g. \cite{A1976,Ten1995}. It is
evident that the exponential Dirichlet series posesses properties
similar to the usual Dirichlet series.

\begin{theorem}
The product of two exponential Dirichlet series is the exponential
Dirichlet series of the binomial convolution of the corresponding
arithmetical functions, i. e.
\begin{equation*}
\widetilde{D}(f,s)\widetilde{D}(g,s)= \widetilde{D}(f\circ g,s).
\label{product}
\end{equation*}
\end{theorem}

{\bf Proof.} According to \eqref{isom},
\begin{equation*}
\widetilde{D}(f\circ g,s)=\widetilde{D}\left(\xi
\left(\frac{f}{\xi}*\frac{g}{\xi}\right),s\right)=
D\left(\frac{f}{\xi}*\frac{g}{\xi},s\right)=
D\left(\frac{f}{\xi},s\right) D\left(\frac{g}{\xi},s\right)
=\widetilde{D}(f,s)\widetilde{D}(g,s).
\end{equation*}

\begin{remark}\upshape
The algebra $(\widetilde{\cal D},+,\cdot, \C)$ of exponential
Dirichlet series is isomorphic to the algebras given above.
\end{remark}

Now, we consider the exponential Dirichlet series of completely
multiplicative functions. Let $\exp(t)=\sum_{k=0}^{\infty} t^k/k!$
be the (formal) exponential power series.

\begin{theorem}\label{th:EDScompl}
If $f$ is completely multiplicative, then
\[
\widetilde{D}(f,s)=\exp \left(\sum_p \frac{f(p)}{p^s}\right).
\]
\end{theorem}

{\bf Proof.} The function $\xi(n)$ is multiplicative. Therefore
using the (formal) Euler product formula,
\begin{eqnarray*}
\widetilde{D}(f,s) &=&\prod_p
\sum_{\nu=0}^{\infty}\frac{f(p^\nu)}{\xi(p^\nu) p^{\nu s}}
= \prod_p\sum_{\nu=0}^{\infty} \frac{f(p)^\nu}{\nu! p^{\nu s}} \\
&=& \prod_p \sum_{\nu=0}^{\infty} \frac1{\nu !} (f(p)p^{-s})^\nu =
\prod_p \exp (f(p)p^{-s}) = \exp  \left(\sum_p
\frac{f(p)}{p^s}\right).
\end{eqnarray*}

This completes the proof.

\bigskip

Let $\zeta_{\P}(s)= \sum_p 1/{p^s}$ denote the prime zeta function
and let $\widetilde{\zeta}(s) = \exp(\zeta_{\P}(s))$ ($\RE s>1$).
Then, by applying Theorem 5.2 and the Glaisher formula
$$\zeta_{\P}(s)=\sum_{n=1}^{\infty} \frac{\mu(n)}n \log
\zeta(ns)\ \ (\RE s>1),$$
we obtain the following consequence.

\begin{corollary}
For $f(n)=n^r$ ($r\in \R$),
\begin{equation*}
\widetilde{D}(f,s) = \widetilde{\zeta}(s-r) =
\prod_{n=1}^{\infty}\zeta(n(s-r))^{\mu(n)/n} \quad (\RE s>r+1).
\end{equation*}
In particular,
\begin{equation*}
\widetilde{\zeta}(s)=\prod_{n=1}^{\infty} \zeta(ns)^{\mu(n)/n} \quad
(\RE s>1).
\end{equation*}
\end{corollary}

\begin{corollary}
For $f(n)=r^{\Omega(n)}$ ($r\in \R$),
\[
\widetilde{D}(r^{\Omega},s) = \exp(r\zeta_{\P}(s))  \quad (\RE s>1).
\]
In particular, for the Liouville function $f(n)=\lambda(n)=
(-1)^{\Omega(n)}$,
\[ \widetilde{D}(\lambda,s)= \exp(-\zeta_{\P}(s)) \quad (\RE s>1).
\]
\end{corollary}

\vskip1mm We now consider a binomial analog of the von Mangoldt
function $\Lambda$. Let
\begin{equation}
\widetilde{\Lambda}(n)=\begin{cases} \log p, \ n=p \text{ prime,} \\
0,  \ \text{otherwise.}
\end{cases}
\end{equation}
Then $\widetilde{\Lambda} \circ I=\log$, since for every $n\ge 1$,
\begin{eqnarray*}
(\widetilde{\Lambda} \circ I)(n) &=&\sum_{d\mid n} \left( \prod_p
\binom{\nu_p(n)}{\nu_p(d)} \right) \widetilde{\Lambda}(d) =\sum_p
\binom{\nu_p(n)}{1} \log p\\
&=& \sum_p \nu_p(n) \log p = \log \prod_p p^{\nu_p(n)}=\log n.
\end{eqnarray*}
Furthermore, $\widetilde{\Lambda}= \lambda \circ \log$ and we obtain
the identities: for any $n$ composite number (i.e. with
$\Omega(n)>1$),
\begin{eqnarray}
\sum_{d\mid n} \left( \prod_p \binom{\nu_p(n)}{\nu_p(d)} \right)
(\log d) (-1)^{\Omega(n/d)} &=& 0,\\
\sum_{d\mid n} \left( \prod_p \binom{\nu_p(n)}{\nu_p(d)} \right)
(-1)^{\Omega(d)}\log d &=& 0.
\end{eqnarray}
We also have
\begin{equation}
\widetilde{D}(\widetilde{\Lambda},s)= \widetilde{D}(\lambda,s)
\widetilde{D}(\log,s)= -
\frac{{\widetilde{\zeta}}'(s)}{\widetilde{\zeta}(s)}, \ \RE s>1.
\end{equation}
Note that for the Chebyshev functions $\theta(x)=\sum_{p\le x} \log
p$ and $\psi(x)=\sum_{p^\nu \le x} \log p$ we have
\begin{equation}
\theta(x)=\sum_{n\le x} \widetilde{\Lambda}(n), \ \
\psi(x)=\sum_{n\le x} \Lambda(n).
\end{equation}


\section{Exponential generating functions}

The generating function or the (formal) power series of $f\in {\cal
A}$ is given by
\begin{equation}
P(f,z)= \sum_{n=1}^{\infty} f(n)z^n. \label{powerseries}
\end{equation}
It is well known that if $f, g\in {\cal A}$, then
\begin{equation}
P(f*g,z)=\sum_{k=1}^{\infty} f(k)P(g,z^k) \label{convoprop}
\end{equation}
formally or assuming that $\sum_{m=1}^{\infty} \sum_{n=1}^{\infty}
f(m)g(n)z^{mn}$ is absolutely convergent.

We define the exponential generating function (or the formal
exponential power series) of $f$ by
\begin{equation}
\widetilde{P}(f,z)= P\left(\frac{f}{\xi},z \right)=
\sum_{n=1}^{\infty} \frac{f(n)}{\xi(n)} z^n. \label{exppowerseries}
\end{equation}
The function $\widetilde{P}(f,z)$ has not hitherto been studied in
the literature but it is analogous to the concept of the exponential
generating function (egf) of a sequence $(a_n)_{n\ge 0}$ given by
\begin{equation}
\widehat{P}(f,z)= \sum_{n=0}^{\infty} \frac{a_n}{n!} z^n,
\label{egf}
\end{equation}
used in combinatorics, see e.g. \cite{S1999}.
Now, $\widetilde{P}(\xi,z)=P(I,z)=\exp(z)$ and we let
\begin{equation}
\Xi(z)=\widetilde{P}(I,z)= \sum_{n=1}^{\infty} \frac1{\xi(n)} z^n
\label{Xi}
\end{equation}
denote the exponential generating function of the constant function
$1$.

For $z=1$,
\[
\sum_{n=1}^{\infty} \frac1{\xi(n)} \ge \sum_p  \frac1{\xi(p)}
=\sum_p 1=\infty,
\]
and for $z=x\in (0, 1)$,
\begin{equation}
0< \Xi(x)= \sum_{n=1}^{\infty} \frac1{\xi(n)} x^n \le
\sum_{n=1}^{\infty} x^n =\frac{x}{1-x};
\end{equation}
hence the radius of convergence of the power series $\Xi(z)$ is
$r=1$.

Also $\Xi(x)=\sum_{n=1}^{\infty} \frac1{\xi(n)}x^n > x+x^2+x^3$ for
any $0<x<1$.

The function $\Xi(z)$ plays here the role of the $\exp$ function.

\begin{theorem}
If $f,g\in {\cal A}$, then
\begin{equation*}
\widetilde{P}(f\circ g,z)= \sum_{k=1}^{\infty} \frac{f(k)}{\xi(k)}
\widetilde{P}\left(g,z^k \right),
\end{equation*}
and the convergence is absolute assuming that $\ds
\sum_{m=1}^{\infty} \sum_{n=1}^{\infty}
\frac{f(m)g(n)}{\xi(m)\xi(n)}z^{mn}$ is absolutely convergent.
\end{theorem}

{\bf Proof.} According to \eqref{isom} and \eqref{convoprop},
\begin{equation*}
\widetilde{P}(f\circ g,z)= \widetilde{P} \left(\xi
\left(\frac{f}{\xi}*\frac{g}{\xi}\right),z \right)=
P\left(\frac{f}{\xi}*\frac{g}{\xi},z\right)= \sum_{k=1}^{\infty}
\frac{f(k)}{\xi(k)}P\left(\frac{g}{\xi},z^k \right)=
\sum_{k=1}^{\infty} \frac{f(k)}{\xi(k)} \widetilde{P}\left(g,z^k
\right).
\end{equation*}
This completes the proof.

\bigskip

 From Corollary 3.1 we see that if $f(n)=r^{\Omega(n)}$
and $g(n)=s^{\Omega(n)}$, then $(f\circ g)(n)=(r+s)^{\Omega(n)}$.

\begin{corollary}
If $f(n)=r^{\Omega(n)}$ and $g(n)=1$, then $(f\circ
g)(n)=(r+1)^{\Omega(n)}$ and
\begin{equation}
\sum_{n=1}^{\infty} \frac{(r+1)^{\Omega(n)}}{\xi(n)} z^n =
\sum_{n=1}^{\infty} \frac{r^{\Omega(n)}}{\xi(n)} \Xi(z^n), \ |z|<1.
\label{Xir}
\end{equation}
In particular, if $f(n)=(-1)^{\Omega(n)}$ and $g(n)=1$, then $f\circ
g=\delta$ and
\begin{equation}
\sum_{n=1}^{\infty} \frac{(-1)^{\Omega(n)}}{\xi(n)} \Xi(z^n)=z, \
|z|<1. \label{Xi1}
\end{equation}
If $f(n)=1$ and $g(n)=1$, then $(f\circ g)(n)=2^{\Omega(n)}$ and
\begin{equation}
\sum_{n=1}^{\infty} \frac{2^{\Omega(n)}}{\xi(n)} z^n =
\sum_{n=1}^{\infty} \frac1{\xi(n)} \Xi(z^n), \ |z|<1. \label{Xi2}
\end{equation}
\end{corollary}


\section{A generalized binomial convolution}

Let $\varphi: \N \times X \to X$, $X\subseteq \C$, be a function
such that, using the notation $\varphi(n,x)=\varphi_n(x)$,

(i) $\varphi_m(\varphi_n(x))=\varphi_{mn}(x)$, $\forall \, m,n\in
\N$, $\forall \,  x\in X$,

(ii) $\varphi_1(x)=x$, $\forall \, x\in X$.

With the aid of the function $\varphi$ we define the following
operation: if $f:\N \to \C$ is an arithmetical function and $\alpha:
X \to \C$ is an arbitrary function, we define $f\boxdot_{\varphi}
\alpha$ by
\begin{equation}
(f \boxdot_{\varphi} \alpha)(x)=\sum_{n=1}^{\infty}
\frac{f(n)}{\xi(n)} \alpha(\varphi_n(x)), \ \forall \, x\in X,
\label{binomoperation}
\end{equation}
assuming that the series is (absolutely) convergent.

If $\varphi_n(x)=x^n$ and $\alpha(x)=x$, then we obtain the
exponential generating function of $f$ given in
\eqref{exppowerseries}. Other special cases considered in this paper
are $\varphi_n(x)=x/n$, $\varphi_n(x)=nx$.

The operation $f\odot_{\varphi} \alpha$ defined by
\begin{equation}
(f \odot_{\varphi} \alpha)(x)=\sum_{n=1}^{\infty} f(n)
\alpha(\varphi_n(x)), \ \forall \, x\in X, \label{convooperation}
\end{equation}
is investigated in the recent paper \cite{BNV2008}, where a detailed
study of \eqref{convooperation} is done, including the problem of
convergence and various applications to M\"obius-type inversion
formulas, even in a more general algebraic context (involving
arithmetical semigroups). A function $\varphi$ satisfying conditions
(i) and (ii) is called an action or a flow on $\N$, the latter term
being used in the theory of dynamical systems.

The following properties of $f \boxdot_{\varphi} \alpha$ are similar
to the properties of $f \odot_{\varphi} \alpha$. In what follows we
write $\boxdot$ instead of $\boxdot_{\varphi}$ for the sake of
brevity.

\begin{theorem}
Let $f,g:\N \to \C, \alpha, \beta: X\to \C$ be arbitrary functions.
Then

1) $f\boxdot (\alpha+\beta)=f \boxdot \alpha + f \boxdot \beta$,

2) $(f+g)\boxdot \alpha= f\boxdot \alpha +g\boxdot \alpha$,

3) $f \boxdot (g \boxdot \alpha)= (f \circ g)\boxdot \alpha$,

4) $\delta \boxdot \alpha  =\alpha$,

\noindent assuming that the appropriate series converge absolutely.
\end{theorem}

{\bf Proof.} Parts 1) and 2) are immediate by the definition. Part
4) is a consequence of (ii), since $(\delta \boxdot \alpha)(x)=
\delta(1)\alpha(\varphi_1(x))=\alpha(x)$, $\forall x\in X$.

Part 3) follows from the similar property $f \odot_{\varphi} (g
\odot_{\varphi} \alpha)= (f*g)\odot_{\varphi} \alpha$, see
\cite[Theorem 1]{BNV2008} and the relation \eqref{isom} between the
Dirichlet convolution and the binomial convolution, but we give here
a direct proof. Using (i),
\begin{eqnarray*}
(f\boxdot (g\boxdot \alpha))(x) &=&\sum_{n=1}^{\infty}
\frac{f(n)}{\xi(n)} (g\boxdot \alpha)(\varphi_n(x))
=\sum_{n=1}^{\infty} \frac{f(n)}{\xi(n)} \sum_{m=1}^{\infty}
\frac{g(m)}{\xi(m)}
\alpha(\varphi_m(\varphi_n(x)))\\
&=&\sum_{n=1}^{\infty} \sum_{m=1}^{\infty}
\frac{f(n)g(m)}{\xi(n)\xi(m)} \alpha(\varphi_{mn}(x)).
\end{eqnarray*}
Assuming that this series is absolutely convergent and grouping its
terms according to the value $nm=k$,
\begin{eqnarray*}
(f\boxdot (g\boxdot \alpha))(x) &=&\sum_{k=1}^{\infty}
\left(\sum_{nm=k} \frac{f(n)g(m)}{\xi(n)\xi(m)} \right)
\alpha(\varphi_k(x)) = \sum_{k=1}^{\infty} \left(\frac{f}{\xi}*
\frac{g}{\xi}\right)(k) \alpha(\varphi_k(x))\\
&=& \sum_{k=1}^{\infty} \frac{(f\circ g)(k)}{\xi(k)}
\alpha(\varphi_k(x))= ((f\circ g)\boxdot \alpha)(x).
\end{eqnarray*}

\begin{theorem}[M\"obius-type inversion]
Let $f\in {\cal A}$ with $f(1)\ne 0$ and let $\alpha, \beta: X\to\C$
be arbitrary functions. Assume that
$$\ds \sum_{n=1}^{\infty}
\sum_{m=1}^{\infty} \frac{f(n)f^{-1\circ}(m)}{\xi(n)\xi(m)}
\beta(\varphi_{mn}(x))$$ is absolutely convergent.
If
\begin{equation}
\alpha(x)=\sum_{n=1}^{\infty} \frac{f(n)}{\xi(n)}
\beta(\varphi_n(x)), \ \forall \, x\in X,  \label{M1} \end{equation}
then
\begin{equation} \beta(x)=\sum_{n=1}^{\infty}
\frac{f^{-1\circ}(n)}{\xi(n)} \alpha(\varphi_n(x)), \ \forall \,
x\in X. \label{M2}
\end{equation}
\end{theorem}

{\bf Proof.} Equation \eqref{M1} can be written as  $\alpha=f\boxdot
\beta$ and therefore using the above theorem we obtain
$f^{-1\circ}\boxdot \alpha =f^{-1\circ}\boxdot (f\boxdot
\beta)=(f^{-1\circ}\circ f)\boxdot \beta = \delta \boxdot \beta=
\beta$. That is, \eqref{M2} holds.

\begin{example}\upshape
If $\alpha(x)=\Xi(x)$, $\beta(x)=x$, $f(n)=1$ and $\varphi_n(x)=x^n$
in \eqref{M1}, then \eqref{M2} becomes \eqref{Xi1}.
\end{example}

We next consider another M\"obius-type inversion. In fact, we
consider the case $\varphi_n(x)=x/n$ and functions $\alpha, \beta:
(0,\infty)\to\C$ such that $\alpha(x)=0$, $\beta(x)=0$ for $x\in
(0,1)$. Here the sums are finite and therefore we need not take care
of convergence.

\begin{theorem}
Let $f\in {\cal A}$ with $f(1)\ne 0$. If
\begin{equation}
\alpha(x)=\sum_{n\le x} \frac{f(n)}{\xi(n)} \beta(x/n), \ \forall \,
x\ge 1,
\end{equation}
then
\begin{equation} \beta(x)=\sum_{n\le x} \frac{f^{-1\circ}(n)}{\xi(n)}
\alpha(x/n), \ \forall \, x\ge 1.
\end{equation}
\end{theorem}

We now consider a M\"obius-type inversion involving only
arithmetical functions. In this case $\varphi_n(x)=nx$ and $x=m\in
\N$. See \cite{LS1980} and \cite{McC1986} for the ``usual'' form.

\begin{theorem}
Let $f,g,h\in {\cal A}$ such that
$$\ds \sum_{k=1}^{\infty}
\sum_{m=1}^{\infty} \frac{h(k)h^{-1\circ}(m)}{\xi(k)\xi(m)} g(kmn)$$
is absolutely convergent for all $n\ge 1$. If
\[
f(n)= \sum_{m=1}^{\infty} \frac{h(m)}{\xi(m)} g(mn), \ \forall \,
n\ge 1,
\]
then
\[
g(n)=\sum_{m=1}^{\infty} \frac{h^{-1\circ}(m)}{\xi(m)} f(mn), \
\forall \, n\ge 1.
\]
\end{theorem}

Theorem 7.4 applies for $h(n)=1$, $h^{-1\circ}(n)=(-1)^{\Omega(n)}$,
also for $h(n)=(-1)^{\Omega(n)}$, $h^{-1\circ}(n)=1$, leading to the
following symmetrical version of the corresponding M\"obius
inversion.

\begin{corollary}
Let $f,g\in {\cal A}$. Then the following two statements are
equivalent:

A) $\ds f(n)= \sum_{m=1}^{\infty} \frac1{\xi(m)} g(mn)$, for every
$n\ge 1$, and $\ds \sum_{n=1}^{\infty} n^{\varepsilon}
|g(n)|<\infty$, for an $\varepsilon >0$.

B) $\ds g(n)= \sum_{m=1}^{\infty} \frac{(-1)^{\Omega(m)}}{\xi(m)}
f(mn)$, for every $n\ge 1$, and $\ds \sum_{n=1}^{\infty}
n^{\varepsilon} |f(n)| <\infty$, for an $\varepsilon >0$.
\end{corollary}

{\bf Proof.} We show that A) implies B). Assume that A) holds. In
order to apply Theorem 7.4 for the appropriate functions we have to
show the series $\sum_{k=1}^{\infty} \sum_{m=1}^{\infty}
g(kmn)/{\xi(k)\xi(m)} $ is absolutely convergent for every $n\ge 1$.
We use that for any $\varepsilon >0$ there is a constant
$C=C(\varepsilon)$ such that $\tau(n):=\sum_{d\mid n} 1 \le
Cn^{\varepsilon}$ for every $n\ge 1$. Thus, for all $K, M\ge 1$, by
grouping the terms according to the value $km=\ell$,
\begin{eqnarray*}
\sum_{k=1}^{K} \sum_{m=1}^{M} \frac1{\xi(k)\xi(m)} |g(kmn)|
&\le&\sum_{\ell=1}^{KM} |g(n\ell)|
\sum_{km=\ell}\frac1{\xi(k)\xi(m)} \le \sum_{\ell=1}^{KM}
|g(n\ell)| \tau(\ell) \\
&\le &  C \sum_{\ell=1}^{KM} |g(n\ell)| \ell^{\varepsilon}=
Cn^{-\varepsilon} \sum_{\ell=1}^{KM} |g(n\ell)|
(n\ell)^{\varepsilon} \\ &\le & Cn^{-\varepsilon}
\sum_{j=1}^{\infty} |g(j)| j^{\varepsilon} < \infty.
\end{eqnarray*}
This shows that the series $\sum_{k=1}^{\infty} \sum_{m=1}^{\infty}
|g(kmn)|/(\xi(k)\xi(m))$ is absolutely convergent for every $n\ge
1$.

Now, applying Theorem 7.4 for $h(n)=1$, $h^{-1\circ}(n)=
(-1)^{\Omega(n)}$ we obtain the first part of B). By similar
arguments as above, for all $N\ge 1$,
\begin{eqnarray*}
\sum_{n=1}^{N} n^{\varepsilon/2} |f(n)| &\le& \sum_{n=1}^{N}
n^{\varepsilon/2} \sum_{m=1}^{\infty} |g(mn)|\le
\sum_{\ell=1}^{\infty} |g(\ell)| \sum_{mn=\ell} n^{\varepsilon/2} \\
&\le& \sum_{\ell=1}^{\infty} |g(\ell)| \ell^{\varepsilon/2}
\tau(\ell) \le C' \sum_{\ell=1}^{\infty} |g(\ell)|
\ell^{\varepsilon}<\infty,
\end{eqnarray*}
where $C'=C(\varepsilon/2)$, proving the second part of B).

\begin{remark}\upshape
A survey of various M\"obius-type functions has been presented in
\cite[Chapter 2]{SC2004}.
\end{remark}

\noindent{\bf Acknowledgement.} We thank Saku Sairanen for providing
a proof of the positive part of Theorem~\ref{th:charS} We also
thank the referee for very careful reading of the manuscript and
useful comments.


\label{lastpage-01}

\end{document}